\documentclass[11pt,a4paper]{article}
\usepackage[colorlinks]{hyperref}

\usepackage[table]{xcolor}
\usepackage{slashed}
\usepackage{xkeyval}
\usepackage{url,amsmath,ifthen,amssymb,latexsym,pstricks,mathrsfs,comment,amsthm,graphicx,tikz,tikz-cd,enumerate,accents,pgffor,cite,wrapfig,multicol,float}

\usepackage{bibspacing}

\usepackage{tkz-euclide}
\usetkzobj{all}

\usepackage{geometry} 

\begin{document}

\newcommand{\nc}{\newcommand}
\nc{\rnc}{\renewcommand}

\let\oldproofname=\proofname
\rnc{\proofname}{\rm\bf{\oldproofname}}

\rnc{\S}{\mathcal S}
\nc{\C}{\mathcal C}
\nc{\D}{\mathcal D}
\nc{\tfnt}{\lfloor\tfrac n2\rfloor}
\nc{\tfmt}{\lfloor\tfrac m2\rfloor}
\nc{\fnt}{\lfloor\frac n2\rfloor}
\nc{\fnf}{\lfloor\frac n4\rfloor}
\nc{\fmt}{\lfloor\frac m2\rfloor}
\nc{\nt}{\tfrac n2}
\nc{\mmt}{\frac m2}
\nc{\mmmt}{\frac m2}
\nc{\ntf}{\lfloor\frac{n+2}4\rfloor}
\nc{\Divm}{\mathbb D(m)}
\nc{\ind}{\chi_\Z}
\nc{\ds}{\displaystyle}
\nc{\pre}{\mid}
\nc{\Mod}[1]{\ (\mathrm{mod}\ #1)}
\rnc{\iff}{\ \Leftrightarrow\ }
\nc{\pfcase}[1]{\medskip \noindent {\bf Case #1.}}
\nc{\eitmc}{\end{itemize}\vspace{-2truemm}}
\rnc{\emptyset}{\varnothing}

\rnc{\P}[2]{P_{#1,#2}}
\nc{\Pmn}{\P{m}{n}}
\nc{\Gmn}{G_{m,n}}
\nc{\Bmn}{B_{m,n}}
\nc{\Tmn}{T_{m,n}}

\nc{\fix}{\operatorname{fix}}
\nc{\Fix}{\operatorname{Fix}}
\nc{\la}{\langle}
\nc{\ra}{\rangle}
\nc{\si}{\sigma}
\nc{\sm}{\setminus}
\nc{\set}[2]{\{ {#1} : {#2} \}} 
\nc{\bigset}[2]{\big\{ {#1} : {#2} \big\}} 
\nc{\pf}{\begin{proof}}
\nc{\epf}{\end{proof}}
\nc{\epfres}{\hfill\qed}
\nc{\epfeq}{\tag*{\qed}}
\nc{\N}{\mathbb N}
\nc{\Z}{\mathbb Z}
\nc{\mt}{\mapsto}
\nc{\id}{\mathrm{id}}
\nc{\COMMA}{,\qquad}
\nc{\COMMa}{,\ \ \ }
\nc{\AND}{\qquad\text{and}\qquad}
\nc{\ANd}{\quad\text{and}\quad}
\nc{\DEFBY}{\qquad\text{defined by}\qquad}
\nc{\BY}{\qquad\text{by}\qquad}
\nc{\ba}{{\bf a}}
\nc{\bb}{{\bf b}}
\nc{\Sum}{\operatorname{\textstyle{\sum}}}

\nc{\bit}{\begin{itemize}}
\nc{\eit}{\end{itemize}}
\nc{\bmc}{\begin{multicols}}
\nc{\emc}{\end{multicols}}
\nc{\itemit}[1]{\item[\emph{(#1)}]}
\nc{\itemnit}[1]{\item[(#1)]}
\nc{\pfitem}[1]{\medskip \noindent (#1).}

\nc{\trianglecircle}[1]{
\tikzstyle{vertex}=[circle,draw=black, fill=white, inner sep = 0.06cm]
\tikzstyle{selectedvertex}=[ultra thick,circle,draw=black, fill=white, inner sep = 0.06cm]
\draw (0,0) circle (2cm);
{\foreach \x in {1,...,9}
{\node[vertex] (\x) at ({108-36*\x}:2) {{\footnotesize $\x$}};}}
\node[vertex] (10) at ({108-36*10}:2) {{\tiny $10$}};
{\foreach \x in {#1}
{\node[selectedvertex] (\x) at ({108-36*\x}:2) {{\footnotesize $\x$}};}}
}

\nc{\trianglecirclereduced}[1]{
\tikzstyle{vertex}=[circle,draw=black, fill=white, inner sep = 0.06cm]
\draw (0,0) circle (2cm);
{\foreach \x in {#1}
{\node[vertex] (\x) at ({108-36*\x}:2) {{\footnotesize $\x$}};}}
}

\nc{\oddcirclediagramnew}[4]{
\tikzstyle{vertex}=[circle,draw=black, fill=white, inner sep = 0.06cm]
\tikzstyle{greenvertex}=[ thick,circle,draw=black, fill=green!20, inner sep = 0.06cm]
\tikzstyle{bluevertex}=[ thick,circle,draw=black, fill=blue!30, inner sep = 0.06cm]
\draw[dashed] (90:2.75)--(270:2.75);
\draw (0,0) circle (2cm);
{\foreach \x/\y in {#2}
{\node[greenvertex] (\x) at ({90 - (360/#1)*(\x)}:2) {{\footnotesize $a_{\y}$}};}}
{\foreach \x/\y in {#3}
{\node[bluevertex] (\x) at ({90 - (360/#1)*(\x)}:2) {{\footnotesize $a_{\y}$}};}}
{\foreach \x/\y in {#4}
{\node[vertex] (\x) at ({90 - (360/#1)*(\x)}:2) {{\footnotesize $a_{\y}$}};}}
\draw[{latex}-{latex}] (260:3) to [bend right=15] (280:3);
\node () at (270:3.3){$\si$};
}

\nc{\evencirclediagramnew}[4]{
\tikzstyle{vertex}=[circle,draw=black, fill=white, inner sep = 0.06cm]
\tikzstyle{greenvertex}=[ thick,circle,draw=black, fill=green!20, inner sep = 0.06cm]
\tikzstyle{bluevertex}=[ thick,circle,draw=black, fill=blue!30, inner sep = 0.06cm]
\draw[dashed] (90:2.75)--(270:2.75);
\draw (0,0) circle (2cm);
{\foreach \x/\y in {#2}
{\node[greenvertex] (\x) at  ({90 - (360/#1)*(\x-.5)}:2) {{\footnotesize $a_{\y}$}};}}
{\foreach \x/\y in {#3}
{\node[bluevertex] (\x) at  ({90 - (360/#1)*(\x-.5)}:2) {{\footnotesize $a_{\y}$}};}}
{\foreach \x/\y in {#4}
{\node[vertex] (\x) at  ({90 - (360/#1)*(\x-.5)}:2) {{\footnotesize $a_{\y}$}};}}
\draw[{latex}-{latex}] (260:3) to [bend right=15] (280:3);
\node () at (270:3.3){$\si$};
}

\nc{\evencirclediagramnewnew}[4]{
\tikzstyle{vertex}=[circle,draw=black, fill=white, inner sep = 0.06cm]
\tikzstyle{greenvertex}=[ thick,circle,draw=black, fill=green!20, inner sep = 0.06cm]
\tikzstyle{bluevertex}=[ thick,circle,draw=black, fill=blue!30, inner sep = 0.06cm]
\draw[dashed] (90:2.75)--(270:2.75);
\draw (0,0) circle (2cm);
{\foreach \x/\y in {#2}
{\node[greenvertex] (\x) at  ({90 - (360/#1)*(\x-1)}:2) {{\footnotesize $a_{\y}$}};}}
{\foreach \x/\y in {#3}
{\node[bluevertex] (\x) at  ({90 - (360/#1)*(\x-1)}:2) {{\footnotesize $a_{\y}$}};}}
{\foreach \x/\y in {#4}
{\node[vertex] (\x) at  ({90 - (360/#1)*(\x-1)}:2) {{\footnotesize $a_{\y}$}};}}
\draw[{latex}-{latex}] (260:3) to [bend right=15] (280:3);
\node () at (270:3.3){$\si$};
}

\numberwithin{equation}{section}

\newtheorem{thm}[equation]{Theorem}
\newtheorem{qu}[equation]{Question}
\newtheorem{lemma}[equation]{Lemma}
\newtheorem{cor}[equation]{Corollary}
\newtheorem{prop}[equation]{Proposition}
\newtheorem*{BL}{Burnside's Lemma}
\newtheorem*{NMP}{The $\boldsymbol{n}$-matchstick Problem}

\theoremstyle{definition}

\newtheorem{rem}[equation]{Remark}
\newtheorem{eg}[equation]{Example}

\title{Integer polygons of given perimeter}
\date{}
\author{James East and Ron Niles}
\maketitle

\vspace{-6mm}

\begin{abstract}
A classical result of Honsberger states that the number of incongruent triangles with integer sides and perimeter $n$ is the nearest integer to $\frac{n^2}{48}$ ($n$ even) or $\frac{(n+3)^2}{48}$ ($n$ odd).  We solve the analogous problem for $m$-gons (for arbitrary but fixed $m\geq3$), and for polygons (with arbitrary number of sides).  We also show that the solution to the latter is asymptotic to $\frac{2^{n-1}}n$, and the former (for fixed $m$) to $\frac{2^{m-1}-m}{2^mm!}n^{m-1}$.



\end{abstract}

\section{Introduction and statement of the main results}\label{sect:intro}

This article concerns \emph{integer polygons}:~i.e., polygons whose side-lengths are all integers.  The study of such polygons dates back at least a few millenia, as the Babylonians and Egyptians were interested in right-angled integer triangles.  A comparatively recent 1904 result attributed to Whitworth and Biddle (see \cite[p.~199]{Dickson1966}) states that there are exactly five (incongruent) integer triangles with perimeter equal to area.  Phelps and Fine \cite{PF1971} showed that there is only one with perimeter equal to twice the area, namely the $(3,4,5)$ triangle.  Subbarao \cite{Subbarao1971} and Marsden \cite{Marsden1974} considered analogous problems for other multiples.  Related to these, the following question appears to have been first asked (and answered) by Jordan, Walch and Wisner \cite{JWW1979} in~1979:

\begin{qu}\label{qu:tri}
How many incongruent integer triangles have perimeter $n$?
\end{qu}

Several formulations of the answer exist \cite{Andrews1979,Honsberger1985,Singmaster1990,KM1998,JWW1979}, and we believe the most elegant is the following:

\begin{thm}[Honsberger \cite{Honsberger1985}]\label{thm:tn}
The number of incongruent integer triangles with perimeter $n$ is $\big[ \tfrac{n^2}{48}\big]$ if $n$ is even, or $\big[ \tfrac{(n+3)^2}{48}\big]$ if $n$ is odd.
\end{thm}

Here, $[x]$ denotes the nearest integer to the real number $x$ (if such exists).  A number of proofs of Theorem~\ref{thm:tn} have been given; see for example \cite{Honsberger1985,JM2000,Hirschhorn2000,Hirschhorn2003}, most of which use recurrence relations and/or generating functions, sometimes ingeniously.

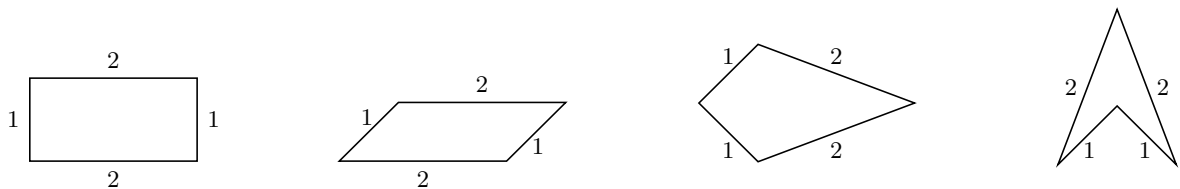
\begin{figure}[ht]%
\begin{center}
\begin{tikzpicture}[scale=1.1]
    \tkzDefPoint(0,0){A}
    \tkzDefPoint(2,0){B}
    \tkzDefPoint(2,1){C}
    \tkzDefPoint(0,1){D}
    \tkzDrawPolygon(A,B,C,D)
    \tkzCalcLength[cm](A,B)\tkzGetLength{ABl}
    \tkzCalcLength[cm](B,C)\tkzGetLength{BCl}
    \tkzCalcLength[cm](C,D)\tkzGetLength{CDl}
    \tkzCalcLength[cm](D,A)\tkzGetLength{DAl}
    \tkzLabelSegment[midway, below](A,B){\footnotesize $\pgfmathprintnumber\ABl$};
    \tkzLabelSegment[midway, right](B,C){\footnotesize $\pgfmathprintnumber\BCl$};
    \tkzLabelSegment[midway, above](C,D){\footnotesize $\pgfmathprintnumber\CDl$};
    \tkzLabelSegment[midway, left](D,A){\footnotesize $\pgfmathprintnumber\DAl$};
\begin{scope}[shift={(3.7,0)}]
    \tkzDefPoint(0,0){A}
    \tkzDefPoint(2,0){B}
    \tkzDefPoint(2+0.70710678118,0.70710678118){C}
    \tkzDefPoint(0.70710678118,0.70710678118){D}
    \tkzDrawPolygon(A,B,C,D)
    \tkzCalcLength[cm](A,B)\tkzGetLength{ABl}
    \tkzCalcLength[cm](B,C)\tkzGetLength{BCl}
    \tkzCalcLength[cm](C,D)\tkzGetLength{CDl}
    \tkzCalcLength[cm](D,A)\tkzGetLength{DAl}
    \tkzLabelSegment[midway, below](A,B){\footnotesize $\pgfmathprintnumber\ABl$};
    \tkzLabelSegment[near start, right](B,C){\footnotesize $\pgfmathprintnumber\BCl$};
    \tkzLabelSegment[midway, above](C,D){\footnotesize $\pgfmathprintnumber\CDl$};
    \tkzLabelSegment[near start, left](D,A){\footnotesize $\pgfmathprintnumber\DAl$};
\end{scope}
\begin{scope}[shift={(8,.7)}]
\begin{scope}[rotate=-45]
    \tkzDefPoint(0,0){A}
    \tkzDefPoint(1,0){B}
    \tkzDefPoint(1.82287565553,1.82287565553){C}
    \tkzDefPoint(0,1){D}
    \tkzDrawPolygon(A,B,C,D)
    \tkzCalcLength[cm](A,B)\tkzGetLength{ABl}
    \tkzCalcLength[cm](B,C)\tkzGetLength{BCl}
    \tkzCalcLength[cm](C,D)\tkzGetLength{CDl}
    \tkzCalcLength[cm](D,A)\tkzGetLength{DAl}
    \tkzLabelSegment[midway, below](A,B){\footnotesize $\pgfmathprintnumber\ABl$};
    \tkzLabelSegment[midway, below](B,C){\footnotesize $\pgfmathprintnumber\BCl$};
    \tkzLabelSegment[midway, above](C,D){\footnotesize $\pgfmathprintnumber\CDl$};
    \tkzLabelSegment[midway, above](D,A){\footnotesize $\pgfmathprintnumber\DAl$};
\end{scope}\end{scope}
\begin{scope}[shift={(13,-.75)}]
\begin{scope}[rotate=45]
    \tkzDefPoint(1,1){A}
    \tkzDefPoint(1,0){B}
    \tkzDefPoint(1.82287565553,1.82287565553){C}
    \tkzDefPoint(0,1){D}
    \tkzDrawPolygon(A,B,C,D)
    \tkzCalcLength[cm](A,B)\tkzGetLength{ABl}
    \tkzCalcLength[cm](B,C)\tkzGetLength{BCl}
    \tkzCalcLength[cm](C,D)\tkzGetLength{CDl}
    \tkzCalcLength[cm](D,A)\tkzGetLength{DAl}
    \tkzLabelSegment[near end, left](A,B){\footnotesize $\pgfmathprintnumber\ABl$};
    \tkzLabelSegment[midway, right](B,C){\footnotesize $\pgfmathprintnumber\BCl$};
    \tkzLabelSegment[midway, left](C,D){\footnotesize $\pgfmathprintnumber\CDl$};
    \tkzLabelSegment[near start, right](D,A){\footnotesize $\pgfmathprintnumber\DAl$};
\end{scope}\end{scope}
\end{tikzpicture}
\end{center}
\vspace{-0.8cm}
\caption{Several quadrilaterals with edge-lengths $1,1,2,2$.}
\label{fig:quad}
\end{figure}

If the word ``triangles'' is simply replaced with ``quadrilaterals'' in Question \ref{qu:tri}, then the answer is not very interesting: for example, there are infinitely many incongruent rhombuses with edges $(1,1,1,1)$.  The same is true for pentagons, hexagons, and so on.  Thus, rather than congruence, we consider a different form of polygon equivalence, defined as follows.  
Let $P$ and $Q$ be $m$-gons for some $m\geq3$, with side lengths $a_1,\ldots,a_m$ and $b_1,\ldots,b_m$, respectively, beginning from any side and reading clockwise or anti-clockwise.  
We say that $P$ and~$Q$ are \emph{equivalent} if we may obtain the $m$-tuple $(b_1,\ldots,b_m)$ from $(a_1,\ldots,a_m)$ by cyclically re-ordering and/or reversing the entries.
In Figure~\ref{fig:quad} for example, the first and second quadrilaterals are equivalent, the third and fourth are equivalent, but the first and third are inequivalent.

Thus, we would like to answer the~following two questions:

\begin{qu}\label{qu:mgon}
How many inequivalent integer $m$-gons have perimeter $n$?
\end{qu}

\begin{qu}\label{qu:poly}
How many inequivalent integer polygons have perimeter $n$?
\end{qu}

As far as we are aware, neither question has been answered previously, apart from the $m=3$ case of Question \ref{qu:mgon} (cf.~Theorem \ref{thm:tn} above); see also \cite{APR2001}, which considers integer $m$-gons up to arbitrary re-orderings of the sides, a different problem for $m\geq4$.
The main purpose of the current article is to answer both Questions~\ref{qu:mgon} and~\ref{qu:poly}, and we now state the results that do so.

If $d$ and $m$ are integers, we write $d\pre m$ to indicate that~$d$ divides $m$: i.e., that $\frac md$ is an integer.  We write~$\lfloor x\rfloor$ for the floor of the real number $x$: i.e., the greatest integer not exceeding $x$.  
A binomial coefficient~$\binom nk$ has its usual meaning if $n,k$ are non-negative integers with $k\leq n$, and is zero otherwise.  Finally, $\varphi$ is Euler's totient function; so for a positive integer~$n$,~$\varphi(n)$ is the size of the set $\set{d\in\{1,\ldots,n\}}{\gcd(d,n)=1}$.  

\begin{thm}\label{thm:pmn}
If $3\leq m\leq n$, then the number $p_{m,n}$ of inequivalent integer $m$-gons with perimeter~$n$~is given by
\[
p_{m,n} = \sum_{d\pre \gcd(m,n)}\frac{\varphi(d)}{2n}\binom{\frac nd}{\frac md} +\frac12 \left(
\binom{\fmt+\lfloor\frac{n-m}2\rfloor}{\fmt}
- \binom{\fnt}{m-1} 
- \binom{\fnf}{\fmt}
- \binom{\ntf}{\mmt}
\right).
\]
\end{thm}

\begin{thm}\label{thm:poly}
If $n\geq3$, then the number $p_n$ of inequivalent integer polygons with perimeter $n$ is given by
\[
p_n = \sum_{d\pre n}\frac{\varphi(d)\cdot2^{\frac nd-1}}n + 2^{\lfloor\frac{n-3}2\rfloor} -
\begin{cases}
3\cdot2^{\lfloor\frac{n-4}4\rfloor} &\text{if $n\equiv0$ or $1\Mod4$}\\[3truemm]
2^{\lfloor\frac{n+2}4\rfloor} &\text{if $n\equiv2$ or $3\Mod4$.} 
\end{cases}
\]
\end{thm}

Although Theorems \ref{thm:pmn} and \ref{thm:poly} are not as striking as Honsberger's $\frac{n^2}{48}$ result, we may nevertheless use these theorems to obtain elegant asymptotic formulae:

\begin{thm}\label{thm:asymptotic}
The number of inequivalent integer polygons with perimeter $n$ is asymptotic to $\frac{2^{n-1}}n$.
\end{thm}

\begin{thm}\label{thm:asymptotic2}
For fixed $m\geq3$, the number of inequivalent integer $m$-gons with perimeter $n$ is asymptotic to $\frac{2^{m-1}-m}{2^mm!}n^{m-1}$.
\end{thm}

From Theorem \ref{thm:pmn}, we may deduce Honsberger's Theorem \ref{thm:tn} as a special case.  As an additional application, we also give an analogous ``nearest integer formula'' for quadrilaterals; again, we are not aware of any previous proof of such a formula.  

\begin{thm}\label{thm:qn}
For any positive integer $n$, the number of inequivalent integer quadrilaterals with perimeter~$n$ is $\big[ \tfrac{n^3-3n^2+20n}{96}\big]$ if $n$ is even, or $\big[ \tfrac{n^3-7n}{96}\big]$ if $n$ is odd.
\end{thm}

Similar formulae could be obtained for pentagons, hexagons, and so on, although these quickly become unweildy.  However, Theorem \ref{thm:asymptotic2} easily leads to the asymptotic formulae:
\[
p_{5,n} \sim \tfrac{11n^4}{3840} \COMMa\ \
p_{6,n} \sim \tfrac{13n^5}{23040} \COMMa\ \
p_{7,n} \sim \tfrac{19n^6}{215040} \COMMa\ \
p_{8,n} \sim \tfrac{n^7}{86016} \COMMa\ \
p_{9,n} \sim \tfrac{247n^8}{185794560} \COMMa\ \
p_{10,n} \sim \tfrac{251n^9}{1857945600}.
\]

The proofs of Theorems \ref{thm:pmn} and \ref{thm:poly} follow the same pattern.  In both cases, we 
(i) identify an action of the dihedral group $\D_n$ (defined below) on a certain set of $n$-tuples, 
(ii) show that the polygons in question are in one-one correspondence with the orbits of the action, and
(iii) enumerate the orbits using Burnside's Lemma (stated below).
We carry out the first two tasks in Section \ref{sect:dihedral}, where we also show how the third reduces to the calculation of certain parameters associated to the elements of $\D_n$.  We calculate these parameters (the bulk of the work) in Section \ref{sect:fix}, and then complete the proofs in Section \ref{sect:proofs}.  Section~\ref{sect:triquad} gives the above-mentioned applications to triangles and quadrilaterals, and Section \ref{sect:values} contains tables of calculated values, a discussion of relevant entries in the Online Encyclopedia of Integer Sequences \cite{OEIS}, and some concluding remarks.

\section{Dihedral group actions}\label{sect:dihedral}

Recall that an action of a group $G$ on a set $X$ is a map $G\times X\to X:(g,x)\mt g\cdot x$, such that $\id_G\cdot x=x$ and $(gh)\cdot x=g\cdot (h\cdot x)$ for all $x\in X$ and $g,h\in G$.  An action induces an equivalence relation on $X$, namely $\bigset{(x,y)\in X\times X}{x=g\cdot y\ (\exists g\in G)}$, the equivalence classes of which are the \emph{orbits} of the action.  The set of orbits is denoted by $X/G$, and the number of orbits is given by Burnside's Lemma (see for example \cite[p246]{PC} for a proof):

\begin{lemma}\label{lem:BS}
If a finite group $G$ acts on a set $X$, then the number of orbits of the action is given by
\[
|X/G| = \frac1{|G|}\sum_{g\in G}\fix_X(g), 
\]
where $\fix_X(g)$ is the cardinality of the set $\Fix_X(g)=\set{x\in X}{g\cdot x=x}$, for $g\in G$. 
\end{lemma}

In what follows, $n$ generally denotes the perimeter of an integer polygon.  Since a polygon has at least three sides, we will always assume that $n\geq3$.  We denote by $\S_n$ the \emph{symmetric group} on the set~$\{1,\ldots,n\}$, which consists of all permutations of this set.
We call a permutation $\si\in\S_n$ a \emph{rotation} if there exists~${q\in\{0,1,\ldots,n-1\}}$ such that $\si(x)\equiv q+x\Mod n$ for all~$x$; when~$q=0$, we obtain the identity map,~$\id_n$.  Similarly,~$\si$ is a \emph{reflection} if there exists~$q$ such that $\si(x)\equiv q-x\Mod n$ for all $x$.  We write~$\D_n$ for the set of all rotations and reflections, and $\C_n$ for the set of all rotations.  These are subgroups of $\S_n$: the \emph{dihedral group} of order $2n$, and the \emph{cyclic group} of order $n$, respectively.  Note that $\D_n$ contains $n$ rotations (including the identity), and $n$ reflections.  When $n$ is odd, all~$n$ reflections fix a single point from $\{1,\ldots,n\}$; when $n$ is even, $\nt$ of the reflections have two fixed points, and~$\nt$ have none.

\begin{figure}[ht]%
\begin{center}
\begin{tikzpicture}[scale=.75]
\trianglecircle{1,3,4,8}
\begin{scope}[shift={(7.5,0)}]
\trianglecirclereduced{1,3,4,8}
\node at ({108-36*6}:2.2) {{\footnotesize $4$}};
\node at ({108-36*2}:2.2) {{\footnotesize $2$}};
\node at ({108-36*3.5}:2.15) {{\footnotesize $1$}};
\node at ({108-36*9.5}:2.2) {{\footnotesize $3$}};
\end{scope}
\begin{scope}[scale=1.131,shift={(11,-.8)}]
\draw (0,0)--(4,0)--(4,1)--(2.255479162,1.978083353)--(0,0)--(4,0);
\node at (2,-.2) {{\footnotesize $4$}};
\node at (3.3,1.6) {{\footnotesize $2$}};
\node at (4.2,.5) {{\footnotesize $1$}};
\node at (1,1.1) {{\footnotesize $3$}};
\end{scope}
\draw[thick,-{latex}] (3,0)--(4.75,0);
\draw[thick,-{latex}] (3+7.5,0)--(4.75+7.5,0);
\begin{scope}[shift={(0,-6)}]
\trianglecircle{3,4,8}
\draw[thick,-{latex}] (3,0)--(4.75,0);
\draw[thick,-{latex}] (3+7.5,0)--(4.75+7.5,0);
\end{scope}
\begin{scope}[shift={(7.5,-6)}]
\trianglecirclereduced{3,4,8}
\node at ({108-36*6}:2.2) {{\footnotesize $4$}};
\node at (90:2.2) {{\footnotesize $5$}};
\node at ({108-36*3.5}:2.15) {{\footnotesize $1$}};
\end{scope}
\begin{scope}[scale=1.131,shift={(11,-6.1)}]
\node at (2.1,.8) {???};
\node[white] at (4.2,.5) {{\footnotesize $1$}};
\end{scope}
\end{tikzpicture}
\end{center}
\vspace{-0.6cm}
\caption{
Top:~the subset $A=\{1,3,4,8\}$ of $\{1,\ldots,10\}$ corresponds to a quadrilateral with side-lengths~$1,2,3,4$.
Bottom:~the subset $B=\{3,4,8\}$ of $\{1,\ldots,10\}$ does not correspond to a polygon.}
\label{fig:circles}
\end{figure}
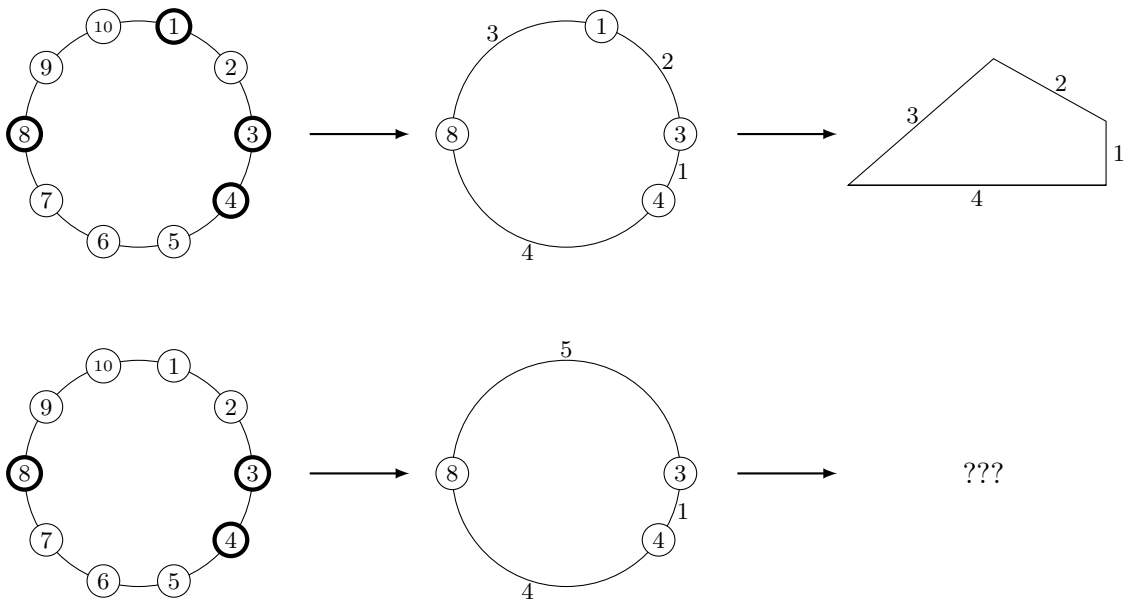

Now consider a circular rope of length $n$ units, with $n$ equally spaced points labelled $1,\ldots,n$, read clockwise.  For any subset~$A$ of $\{1,\ldots,n\}$, we may attempt to create a polygon out of the rope, with corners at the points from~$A$, by pulling the strings taut between these selected points.  With $n=10$, for example, this is possible for $A=\{1,3,4,8\}$, but not for $B=\{3,4,8\}$; see Figure \ref{fig:circles}.  Obviously, to create a polygon for such a subset $A$ of $\{1,\ldots,n\}$, we would need $|A|\geq3$, but there is an additional restriction coming from the fact that the sides of a polygon must all be less than half the perimeter.  Namely, if $A=\{x_1,\ldots,x_m\}$, where $m\geq3$ and $x_1<\cdots<x_m$, then we require
\[
\max\{x_2-x_1,x_3-x_2,\ldots,x_m-x_{m-1},n+x_1-x_m\} < \tfrac n2
\]
to ensure that the length of any side is less than the combined lengths of the other $m-1$ sides.

The above ideas are easier to work with when formulated in slightly different terms.
We define the set $T=\{0,1\}$, and denote by $T_n$ the set of all $n$-tuples over $T$.  We identify a subset $A$ of $\{1,\ldots,n\}$ with an $n$-tuple $\ba=(a_1,\ldots,a_n)\in T_n$ in the usual way: we have $a_x=1 \iff x\in A$.  Thus, given $\ba\in T_n$, we may attempt to form a polygon of perimeter $n$ as above.  It turns out that the $n$-tuples corresponding to polygons in this way may be described very easily.  

By a \emph{block of 0's of length $l$} of an $n$-tuple $\ba\in T_n$, we mean a sequence of $l$ consecutive entries of~$\ba$ (possibly ``wrapping around $n$''), all of which are 0, that is not contained in any larger such sequence of consecutive~0's.
Thus, for example, the subsets $A$ and $B$ of $\{1,\ldots,10\}$ defined above (cf.~Figure \ref{fig:circles}) correspond to the~10-tuples~$\ba = (1,0,1,1,0,0,0,1,0,0)$ and~$\bb = (0,0,1,1,0,0,0,1,0,0)$; here, $\ba$ has three blocks of 0's, of lengths~1,~2 and~3, while $\bb$ has blocks of lengths 3 and 4; one block of $\bb$ wraps around.  

Writing $k=\fnt$, so that $n=2k$ or $2k+1$, we will call a block of 0's of an $n$-tuple $\ba\in T_n$ \emph{bad} if its length is at least $k-1$ ($n$ even) or~$k$ ($n$ odd).  
If an $n$-tuple $\ba\in T_n$ had a bad block of 0's, then it would be impossible to form a polygon from $\ba$, as the bad block would lead to an edge of length at least $\nt$.  Conversely, if $\ba\in T_n$ has no bad blocks, then we can form a polygon from $\ba$; for this, note that having no bad blocks forces $\ba$ to have at least three 1's, since $n\geq3$.  Accordingly, we call an $n$-tuple \emph{bad} if it has at least one bad block, or \emph{good} if it has none, and we write
\[
G_n = \set{\ba\in T_n}{\ba \text{ is good}} \AND B_n = \set{\ba\in T_n}{\ba \text{ is bad}}.
\]
So $G_n$ consists of the $n$-tuples that correspond to polygons in the manner described above.

For $\ba=(a_1,\ldots,a_n)\in T_n$, we denote by $\Sum\ba$ the integer $a_1+\cdots+a_n$, which is just the number of 1's in $\ba$.  For $0\leq m\leq n$, we write
\[
\Tmn=\bigset{\ba\in T_n}{\Sum\ba=m}
\]
for the set of all $n$-tuples with exactly $m$ 1's, and we also write
\[
\Gmn = G_n \cap \Tmn \AND \Bmn = B_n\cap \Tmn
\]
for the sets of all good and bad such $n$-tuples, respectively.  So $\Gmn$ consists of the $n$-tuples that correspond to $m$-gons in the manner described above.  Note that $\Gmn=\emptyset$ if $m\leq2$.

There is a natural action of the dihedral group $\D_n$ on $T_n$ given by permuting the coordinates:
\begin{equation}\label{eq:action}
\si\cdot(a_1,\ldots,a_n) = (a_{\si^{-1}(1)},\ldots,a_{\si^{-1}(n)}) \qquad\text{for $(a_1,\ldots,a_n)\in T_n$ and $\si\in\D_n$.}
\end{equation}
It is easy to see that for any $\ba\in T_n$ and $\si\in\D_n$, we have
\bit\bmc3
\itemnit{i} $\Sum\ba=\Sum(\si\cdot\ba)$,
\itemnit{ii} $\ba\in G_n \iff \si\cdot\ba\in G_n$,
\itemnit{iii} $\ba\in B_n \iff \si\cdot\ba\in B_n$.
\emc\eit
It follows from (i) that $\Tmn$ is closed under the action of $\D_n$; by~(ii) and (iii), so too are the sets $G_n$ and~$B_n$; by combinations of these facts, it follows that the sets~$\Gmn$ and~$\Bmn$ are closed as well.  Thus, \eqref{eq:action} defines actions of $\D_n$ on all of the sets $T_n$, $G_n$, $B_n$, $\Tmn$, $\Gmn$ and $\Bmn$.

Clearly two good $n$-tuples are in the same orbit of the action if and only if they determine equivalent polygons.  Together with Lemma \ref{lem:BS}, and writing $p_n$ (respectively, $p_{m,n}$) for the number of inequivalent integer polygons (respectively, $m$-gons) of perimeter $n$, it follows that
\begin{equation}\label{eq:BS}
p_n = |G_n/\D_n| = \frac1{2n}\sum_{\si\in\D_n}\fix_{G_n}(\si) \AND p_{m,n} = |\Gmn/\D_n| = \frac1{2n}\sum_{\si\in\D_n}\fix_{\Gmn}(\si).
\end{equation}
Also, for any $\si\in\D_n$, $\Fix_{T_n}(\si) = \Fix_{G_n}(\si) \cup \Fix_{B_n}(\si)$, with a similar statement for $\Fix_{\Tmn}(\si)$.  Since $G_n$ and $B_n$ are disjoint, it follows that 
\begin{equation}\label{eq:fix}
\fix_{G_n}(\si) = \fix_{T_n}(\si) - \fix_{B_n}(\si) \AND \fix_{\Gmn}(\si) = \fix_{\Tmn}(\si) - \fix_{\Bmn}(\si)  \qquad\text{for any $\si\in\D_n$.}
\end{equation}
Equations \eqref{eq:BS} and \eqref{eq:fix} form the basis of our calculation of $p_n$ and $p_{m,n}$.  We use \eqref{eq:fix} to calculate the values of $\fix_{G_n}(\si)$ and $\fix_{\Gmn}(\si)$ in Section \ref{sect:fix}, before showing in Section \ref{sect:proofs} that these, together with \eqref{eq:BS}, yield the formulae for $p_n$ and $p_{m,n}$ stated in Theorems \ref{thm:pmn} and \ref{thm:poly}.

\section{Fix sets}\label{sect:fix}

The previous section reduced the enumeration of integer polygons of given perimeter to the calculation of sizes of fix sets under the action \eqref{eq:action}.  We perform these calculations in the current section.  

Note that each result of this section gives values for $\fix_{G_n}(\si)$ and $\fix_{\Gmn}(\si)$ for various elements $\si\in\D_n$.  In principle, it would be possible to derive the former from the latter by summing over $m$.  However, this is easier said than done in most cases; in any event, both values may be calculated with essentially the same argument, so this is the approach we take.

In all the proofs that follow, if $\ba\in T_r$ for some $r$, we generally assume that $\ba=(a_1,\ldots,a_r)$.  We also remind the reader of the convention regarding binomial coefficients being zero if the arguments fall outside of the usual ranges.  We begin with the identity element.

\begin{lemma}\label{lem:poly1}
If $n\geq3$, then
\bit
\itemit{i}
$
\fix_{G_n}(\id_n) = 2^n-1 - n\cdot2^{\lfloor\frac n2\rfloor} + \begin{cases}
\nt &\text{if $n$ is even}\\
0 &\text{if $n$ is odd,}
\end{cases}
$
\itemit{ii} 
$\ds\fix_{\Gmn}(\id_n) = \binom nm - n\binom{\fnt}{m-1}$ for any $3\leq m\leq n$.
\eit
\end{lemma}

\pf
We write $k=\fnt$ throughout the proof, so that $n=2k$ or $2k+1$.

\pfitem{i}
Clearly $\fix_{G_n}(\id_n)=|G_n|=|T_n|-|B_n|$.  Since $|T_n|=2^n$, it therefore suffices to show that 
\[
|B_n| = 
1+n\cdot2^k +
\begin{cases}
-k &\text{if $n=2k$ is even} \\
0 &\text{if $n=2k+1$ is odd.}
\end{cases}
\]
Certainly $(0,\ldots,0)\in B_n$, so it remains to count the non-zero $n$-tuples from $B_n$.
By definition, every $n$-tuple from $B_n$ has at least one bad block of 0's.  In fact, it is only possible for an $n$-tuple to have two bad blocks if $n$ is even, in which case there are exactly $\nt$ such $n$-tuples (these have two 1's equally spaced around the circle).
For $i\in\{1,\ldots,n\}$, write~$\pi_i$ for the number of non-zero $n$-tuples $\ba$ from $T_n$ that have a bad block beginning at position $i$: i.e., reading subscripts modulo $n$,
\[
a_{i-1}=1 \AND a_i=a_{i+1}=\cdots=a_{i+l-1}=0, \qquad\text{where $l=k-1$ ($n$ even) or $l=k$ ($n$ odd).}
\]
Since the remaining $n-(l+1)=k$ entries $a_{i+l},\ldots,a_{i-2}$ of such an $n$-tuple may be chosen arbitrarily from~$T$, it follows that $\pi_i=2^k$ for all $i$.
Then $\pi_1+\cdots+\pi_n=n\cdot2^k$ counts the non-zero $n$-tuples with one bad block once, but double-counts the non-zero $n$-tuples with two bad blocks.  We have already noted that there are none of the latter if $n$ is odd, or $\nt=k$ of them if $n$ is even.  The result quickly follows.

\pfitem{ii} This time, $\fix_{\Gmn}(\id_n)=|\Gmn|=|\Tmn|-|\Bmn|=\binom nm-|\Bmn|$.  Every element of~$\Bmn$ has exactly one bad block of 0's (since $m\geq3$).  As in the proof of (i), one may show that there are $\binom k{m-1}=\binom{\fnt}{m-1}$ elements of~$\Bmn$ with a bad block beginning at position~$i$, and the result quickly follows.
\epf

Next we consider the non-trivial rotations.  Recall that the \emph{order} of an element $\si\in\D_n$ is the least positive integer $d$ such that $\si^d=\id_n$; the order of a rotation is a divisor of $n$ and, conversely, any such divisor can occur as the order of a rotation.

\begin{lemma}\label{lem:poly2}
If $n\geq3$, and if $\si\in\D_n$ is a rotation of order $d$, where $1\not=d\pre n$, then
\bit
\itemit{i}
$\fix_{G_n}(\si) = 2^{\frac nd}-1+
\begin{cases}
-\frac n2 &\text{if $d=2$}\\
0 &\text{if $d\geq3$,}
\end{cases}
$
\itemit{ii} $\ds\fix_{\Gmn}(\si) = \binom{\frac nd}{\frac md}$ for any $3\leq m\leq n$.
\eit
\end{lemma}

\pf
Let $e=\frac nd$.  Since all rotations of order $d$ are powers of each other, each such rotation fixes the same $n$-tuples from $G_n$ and $\Gmn$.  Thus, we may assume that $\si(x)\equiv x+e\Mod n$ for all~$x$~(cf.~Figure \ref{fig:circles0}).  For use in both parts of the proof, we define a map $f:T_e\to T_n$ by ${f(a_1,\ldots,a_e)=(a_1,\ldots,a_e,a_1,\ldots,a_e,\ldots,a_1,\ldots,a_e)}$.

\pfitem{i}
It is easy to see that the set $\Fix_{T_n}(\si)$ is precisely the image of $f$ (cf.~Figure \ref{fig:circles0}).  Since $f$ is clearly injective, it immediately follows that $\fix_{T_n}(\si)=|T_e|=2^e=2^{\frac nd}$.  By~\eqref{eq:fix}, it remains to show that 
\[
\fix_{B_n}(\si) = 1+
\begin{cases}
\frac n2 &\text{if $d=2$}\\
0 &\text{if $d\geq3$.}
\end{cases}
\]
Clearly $(0,\ldots,0)\in\Fix_{B_n}(\si)$.  If $\ba\in T_e$ is non-zero, then the longest block of 0's in $f(\ba)$ has length at most~$e-1=\frac nd-1$ (cf.~Figure \ref{fig:circles0}).  If $d>2$, then
\[
\tfrac nd-1 < \tfrac n2-1 = \begin{cases}
k-1 &\text{if $n=2k$ is even}\\
k-\tfrac12 &\text{if $n=2k+1$ is odd,}
\end{cases}
\]
and so $f(\ba)$ belongs to $G_n$ for any non-zero $\ba\in T_e$.  It follows that $\fix_{B_n}(\si) = 1$ if $d>2$.  If $d=2$, then $n=2k$ must be even; it is possible for $\ba\in T_e=T_k$ to be non-zero, but to have $f(\ba)\in B_n$; this occurs when~$\ba$ has exactly one non-zero entry (cf.~Figure \ref{fig:circles0}, right).  It follows that $\fix_{B_n}(\si) = 1+k = 1+\nt$ if $d=2$.  

\pfitem{ii} With the above notation, $\Fix_{\Tmn}=\set{f(\ba)}{\ba\in T_{\frac md,\frac nd}}$, and so $\fix_{\Tmn}(\si)=|T_{\frac md,\frac nd}|=\binom{\frac nd}{\frac md}$.
In the proof of (i), we showed that any non-zero elements of $\Fix_{B_n}(\si)$ have exactly two 1's, and hence do not belong to~$\Bmn$ (since $m\geq3$).  It follows that $\fix_{\Bmn}(\si)=0$.  Because of \eqref{eq:fix}, this completes the proof.
\epf

\nc{\rotationdiagramc}[2]{
\tikzstyle{firstvertex}=[ thick,circle,draw=black, fill=blue!41, inner sep = 0.06cm]
\tikzstyle{secondvertex}=[ thick,circle,draw=black, fill=blue!34, inner sep = 0.06cm]
\tikzstyle{thirdvertex}=[ thick,circle,draw=black, fill=blue!27, inner sep = 0.06cm]
\tikzstyle{fourthvertex}=[ thick,circle,draw=black, fill=blue!20, inner sep = 0.06cm]
\tikzstyle{fifthvertex}=[ thick,circle,draw=black, fill=blue!13, inner sep = 0.06cm]
\tikzstyle{sixthvertex}=[ thick,circle,draw=black, fill=blue!06, inner sep = 0.06cm]
\draw (0,0) circle (2cm);
{\foreach \x/\y/\z in {#2}
{\node[\z] (\x) at  ({90 - (360/#1)*(\x-.5)}:2) {{\footnotesize $a_{\y}$}};}}
}

\begin{figure}[t]%
\begin{center}
\begin{tikzpicture}[scale=1]
\begin{scope}
\rotationdiagramc{12}{
1/1/thirdvertex,2/2/fourthvertex,3/3/fifthvertex,4/4/sixthvertex,
5/1/thirdvertex,6/2/fourthvertex,7/3/fifthvertex,8/4/sixthvertex,
9/1/thirdvertex,10/2/fourthvertex,11/3/fifthvertex,12/4/sixthvertex}
\draw[dashed] (90:2.75)--(0,0)--(330:2.75) (0,0)--(210:2.75);
\draw[-{latex}] (90:.75) arc (90:-30:.75);
\node () at (30:1){$\si$};
\end{scope}
\begin{scope}[shift={(8,0)}]
\draw[white] (90:2.75)--(0,0)--(330:2.75) (0,0)--(210:2.75);
\rotationdiagramc{12}{
1/1/firstvertex,2/2/secondvertex,3/3/thirdvertex,4/4/fourthvertex,5/5/fifthvertex,6/6/sixthvertex,
7/1/firstvertex,8/2/secondvertex,9/3/thirdvertex,10/4/fourthvertex,11/5/fifthvertex,12/6/sixthvertex}
\draw[dashed] (90:2.75)--(0,0)--(270:2.75);
\draw[-{latex}] (90:.75) arc (90:-90:.75);
\node () at (0:1){$\si$};
\end{scope}
\end{tikzpicture}
\end{center}
\vspace{-0.6cm}
\caption{The rotation $\si\in\D_n$ and the $n$-tuple $f(\ba)=(a_1,\ldots,a_e,a_1,\ldots,a_e,\ldots,a_1,\ldots,a_e)\in\Fix_{T_n}(\si)$, from the proof of Lemma~\ref{lem:poly2}, in the case $n=12$ and $d=3$ (left) and $d=2$ (right).}
\label{fig:circles0}
\end{figure}
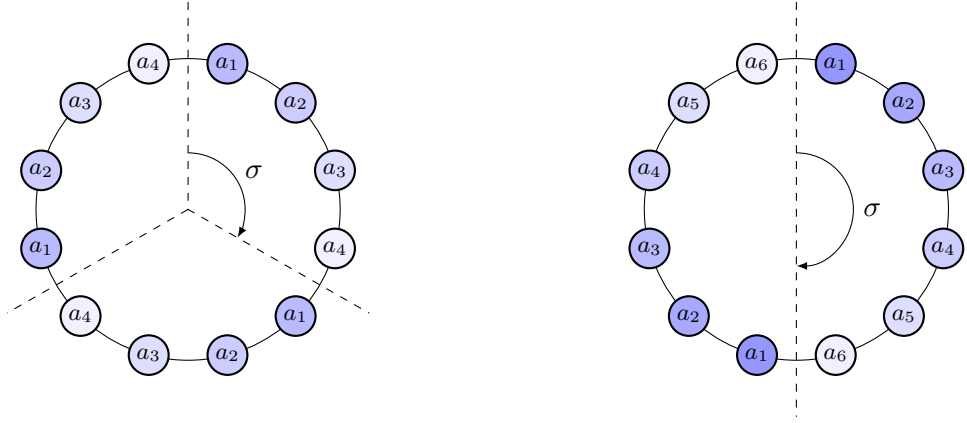

For the reflections, we need to consider separate cases according to the parity of $n$, and according to the number of fixed points in the case of even $n$.  Although the details of the proofs vary, the strategy is essentially the same in each case; we show that an $n$-tuple fixed by a reflection is uniquely determined by (roughly) half of its entries, and identify the properties of these entries that separate the good $n$-tuples from the bad.
We begin with the case of odd $n$.  Recall that here there are $n$ reflections, each of which fixes a single point of $\{1,\ldots,n\}$.

\begin{lemma}\label{lem:poly3}
If $n\geq3$ is odd, and if $\si\in\D_n$ is a reflection, then
\bit
\itemit{i} 
$\ds
\fix_{G_n}(\si) = \begin{cases}
\ds2^{\frac{n+1}2}-3\cdot2^{\frac{n-1}4}+1 &\text{if $n\equiv1\Mod4$}\\
\ds2^{\frac{n+1}2}-2^{\frac{n+5}4}+1 &\text{if $n\equiv3\Mod4$,}
\end{cases}
$
\itemit{ii} $\ds\fix_{\Gmn}(\si) = 
\binom{\fnt}{\fmt} - \binom{\fnf}{\fmt}
- \binom{\ntf}{\frac m2}$ for any $3\leq m\leq n$.
\eit
\end{lemma}

\pf
Write $n=2k+1$, and let $l=\lfloor\frac k2\rfloor$, so that $k=2l$ or $2l+1$.  All reflections fix the same number of $n$-tuples from $G_n$ and $\Gmn$ (and from $T_n$, $B_n$, etc.), so we may assume that the fixed point of $\si$ is $n$, in which case $\si(x)\equiv n-x\Mod n$ for all~$x$ (cf.~Figure \ref{fig:circles2}).  
This time, we define a map $f:T_{k+1}\to T_n$ by~$f(a_1,\ldots,a_{k+1})=(a_1,\ldots,a_k,a_k,\ldots,a_1,a_{k+1})$.  

\pfitem{i}
As in the previous proof, $f$ is injective and its image is $\Fix_{T_n}(\si)$.  Thus, $\fix_{T_n}(\si) = 2^{k+1} = 2^{\frac{n+1}2}$.  By~\eqref{eq:fix}, it remains to show that
\begin{equation}\label{eq:poly3}
\fix_{B_n}(\si) = \begin{cases}
\ds3\cdot2^{\frac{n-1}4}-1 &\text{if $k=2l$}\\
\ds2^{\frac{n+5}4}-1 &\text{if $k=2l+1$.}
\end{cases}
\end{equation}
To do so, consider some $n$-tuple $f(\ba)\in\Fix_{T_n}(\si)$, where $\ba\in T_{k+1}$.  Then $f(\ba)\in B_n$ if and only if at least one of the following holds (cf.~Figure \ref{fig:circles2}):
\begin{itemize}\begin{multicols}{2}
\item[(a)] $(a_1,\ldots,a_l,a_{k+1})=(0,\ldots,0)$, or
\item[(b)] $(a_{l+1},\ldots,a_k)=(0,\ldots,0)$.
\end{multicols}\eitmc
There are $2^{k-l}$ $(k+1)$-tuples~$\ba$ satisfying~(a), $2^{l+1}$ satisfying (b), and one satisfying both (a) and (b).
It follows that $\fix_{B_n}(\si) = 2^{k-l}+2^{l+1}-1$; this reduces to the expressions stated in \eqref{eq:poly3}, by checking separate cases for $k=2l$ or $2l+1$.

\pfitem{ii} Now consider an element $f(\ba)\in\Fix_{\Tmn}(\si)$, where $\ba\in T_{k+1}$.  By the form of $f(\ba)$, and since $n$ is the unique point of $\{1,\ldots,n\}$ fixed by $\si$, we must have
\[
a_{k+1}=\begin{cases}
0 &\text{if $m$ is even}\\
1 &\text{if $m$ is odd.}
\end{cases}
\]
The remaining entries $a_1,\ldots,a_k$ of $\ba$ may be chosen arbitrarily from $T$, as long as $\fmt$ of them are 1's.  Thus, $\fix_{\Tmn}(\si)=\binom k{\fmt}=\binom{\fnt}{\fmt}$.  By \eqref{eq:fix}, it remains to show that $\fix_{\Bmn}(\si)=\binom{\fnf}{\fmt}
+ \binom{\ntf}{\frac m2}$.  
Now,~$f(\ba)$ belongs to $B_n$ if and only if one of (a) or (b) holds, as enumerated in the proof of (i).  Condition (b) holds if and only if $\fmt$ of the entries $a_1,\ldots,a_l$ are 1's (recall that $a_{k+1}$ is already fixed), so there are $\binom l{\fmt}=\binom{\fnf}{\fmt}$ $(k+1)$-tuples~$\ba$, with $f(\ba)\in\Bmn$, satisfying (b).  For $\ba$ to satisfy (a), $m$ must be even (since $a_{k+1}=1$ if~$m$ is odd).  In this case, $\mmt$ of $a_{l+1},\ldots,a_k$ must be 1's; it follows that there are $\binom{k-l}{\mmt}=\binom{\ntf}{\mmmt}$ $(k+1)$-tuples~$\ba$ satisfying (a) when $m$ is even; this is also true when $m$ is odd, since then $\mmt$ is not an integer.  Since no $\ba$ can satisfy both (a) and (b), as $m\geq3$, it follows that $\fix_{\Bmn}(\si)=\binom{\fnf}{\fmt} + \binom{\ntf}{\frac m2}$, as required.
\epf

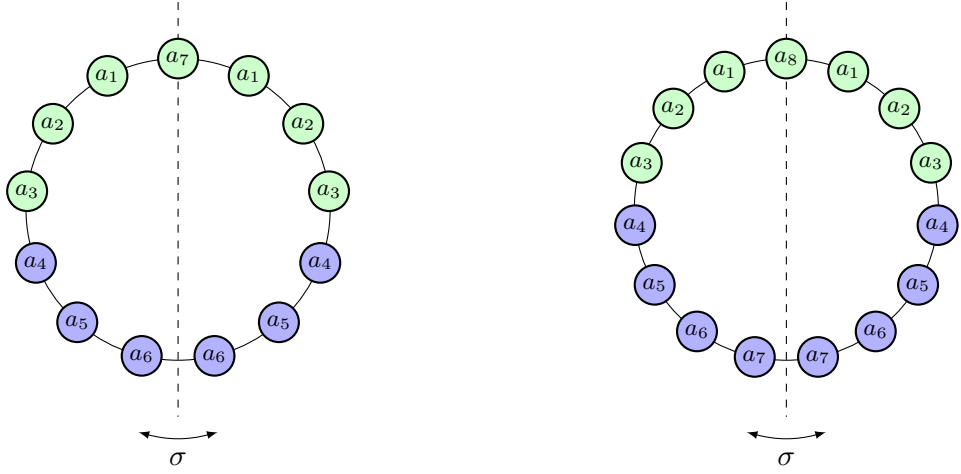
\begin{figure}[t]%
\begin{center}
\begin{tikzpicture}[scale=1]
\begin{scope}
\oddcirclediagramnew{13}{0/7,1/1,2/2,3/3,10/3,11/2,12/1}{4/4,5/5,6/6,7/6,8/5,9/4}{}
\end{scope}
\begin{scope}[shift={(8,0)}]
\oddcirclediagramnew{15}{0/8,1/1,2/2,3/3,12/3,13/2,14/1}{4/4,5/5,6/6,7/7,8/7,9/6,10/5,11/4}{}
\end{scope}
\end{tikzpicture}
\end{center}
\vspace{-0.6cm}
\caption{The reflection $\si\in\D_n$ and the $n$-tuple $f(\ba)=(a_1,\ldots,a_k,a_k,\ldots,a_1,a_{k+1})\in\Fix_{T_n}(\si)$, from the proof of Lemma \ref{lem:poly3}, in the cases $n=13$ (left) and $n=15$ (right).  If condition (a) or (b) from the proof holds, then the green (light) or blue (dark) vertices, respectively, yield a bad block of 0's in $f(\ba)$.  If neither~(a) nor (b) holds, then~$f(\ba)$ has no bad blocks.}
\label{fig:circles2}
\end{figure}

Recall that there are two different kinds of reflections when $n$ is even: $\nt$ fixing no points of $\{1,\ldots,n\}$, and $\nt$ fixing two.  We consider these separately.

\begin{lemma}\label{lem:poly4}
If $n\geq4$ is even, and if $\si\in\D_n$ is a reflection with no fixed points, then
\bit
\itemit{i} $
\fix_{G_n}(\si) = \begin{cases}
\ds 2^{\frac n2} - 2^{\frac{n+4}4}+1 &\text{if $n\equiv0\Mod4$}\\
\ds 2^{\frac n2} - 2^{\frac{n+6}4}+2 &\text{if $n\equiv2\Mod4$,}
\end{cases}
$
\itemit{ii} $\ds\fix_{\Gmn}(\si) = \binom{\nt}{\mmt} -2\binom{\ntf}{\mmt}$ for any $3\leq m\leq n$.
\eit
\end{lemma}

\pf
Write $n=2k$, and let $l=\lfloor\frac k2\rfloor$.  This time we may assume that $\si(x)\equiv n+1-x\Mod n$ for all $x$ (cf.~Figure \ref{fig:circles3}), and we define a map $f:T_k\to T_n$ by $f(a_1,\ldots,a_k)=(a_1,\ldots,a_k,a_k,\ldots,a_1)\in T_n$.  

\pfitem{i}
Again $f$ is injective and has image $\Fix_{T_n}(\si)$, so that $\fix_{T_n}(\si)=2^k=2^{\frac n2}$.  Thus, again by \eqref{eq:fix}, it remains to show that
\begin{equation}\label{eq:poly4}
\fix_{B_n}(\si) = \begin{cases}
\ds2^{\frac{n+4}4}-1 &\text{if $k=2l$}\\
\ds2^{\frac{n+6}4}-2 &\text{if $k=2l+1$.}
\end{cases}
\end{equation}
\newpage\noindent
To do so, consider some $n$-tuple $f(\ba)\in\Fix_{T_n}(\si)$, where $\ba\in T_k$.  Then $f(\ba)\in B_n$ if and only if at least one of the following holds (cf.~Figure \ref{fig:circles3}):
\begin{itemize}\begin{multicols}{2}
\item[(a)] $(a_1,\ldots,a_l)=(0,\ldots,0)$, or
\item[(b)] $(a_{k-l+1},\ldots,a_k)=(0,\ldots,0)$.
\end{multicols}\eitmc
There are $2^{k-l}$ $k$-tuples~$\ba$ satisfying~(a), and also $2^{k-l}$ satisfying (b).  Only one $\ba$ satisfies both (a) and (b) if~$k=2l$, but there are two such $\ba$ if $k=2l+1$ (compare the left and right diagrams in Figure \ref{fig:circles3}).  It follows that $\fix_{B_n}(\si)=2\cdot2^{k-l}-(1\text{ or }2)$, as appropriate, which reduces to \eqref{eq:poly4}.

\pfitem{ii}  Now, $\Fix_{\Tmn}(\si)$ is empty if $m$ is odd (cf.~Figure~\ref{fig:circles3}), in which case the stated formula for $\fix_{\Gmn}(\si)$ holds, as $\mmt$ is not an integer.  From now on, we assume $m=2h$ is even.  Clearly, $\fix_{\Tmn}(\si)=\binom kh$.  If $\ba\in T_k$ is such that $f(\ba)$ belongs to $\Bmn$, then exactly one of (a) or (b) holds, as above.  There are $\binom{k-l}h$ $k$-tuples~$\ba$ satisfying (a), and the same number satisfying (b).  Thus, by \eqref{eq:fix},
\[
\fix_{\Gmn}(\si) = \fix_{\Tmn}(\si) - \fix_{\Bmn}(\si) = \binom kh - 2\binom{k-l}h = \binom{\nt}{\mmt} - 2\binom{\ntf}{\mmt}. \qedhere
\]
\epf

\begin{figure}[t]%
\begin{center}
\begin{tikzpicture}[scale=1]
\begin{scope}
\evencirclediagramnew{12}{1/1,2/2,3/3,10/3,11/2,12/1}{4/4,5/5,6/6,7/6,8/5,9/4}{}
\end{scope}
\begin{scope}[shift={(8,0)}]
\evencirclediagramnew{14}{1/1,2/2,3/3,12/3,13/2,14/1}{5/5,6/6,7/7,8/7,9/6,10/5}{4/4,11/4}
\end{scope}
\end{tikzpicture}
\end{center}
\vspace{-0.6cm}
\caption{The reflection $\si\in\D_n$ and the $n$-tuple $f(\ba)=(a_1,\ldots,a_k,a_k,\ldots,a_1)\in\Fix_{T_n}(\si)$, from the proof of Lemma~\ref{lem:poly4}, in the cases $n=12$ (left) and $n=14$ (right).  If condition (a) or (b) from the proof is satisfied, then the green (light) or blue (dark) vertices, respectively, yield a bad block of 0's in $f(\ba)$.  If neither (a) nor (b) holds, then~$f(\ba)$ has no bad blocks.}
\label{fig:circles3}
\end{figure}
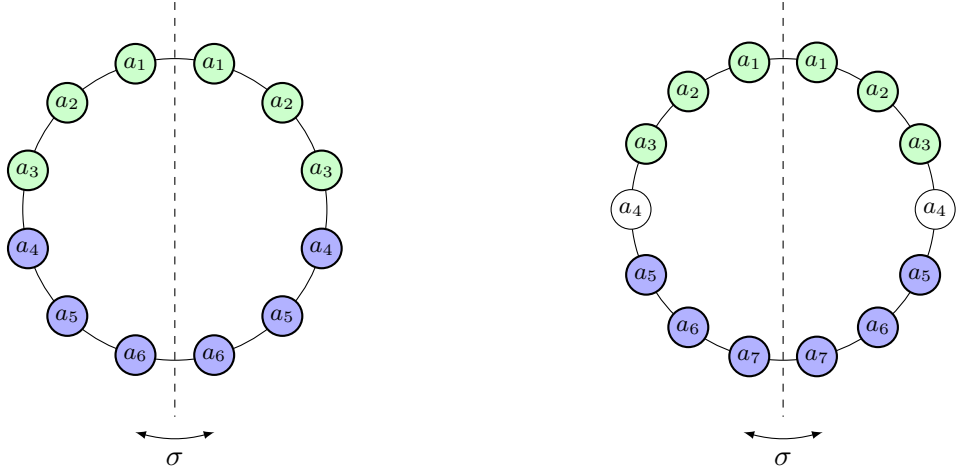

For the final lemma, we must consider separate cases for $\fix_{\Gmn}(\si)$ according to the parity of $m$.

\begin{lemma}\label{lem:poly5}
If $n\geq4$ is even, and if $\si\in\D_n$ is a reflection with two fixed points, then
\bit
\itemit{i} $\fix_{G_n}(\si) = 
\begin{cases}
\ds 2^{\frac{n+2}2} - 2^{\frac{n+8}4}+1 &\text{if $n\equiv0\Mod4$}\\
\ds 2^{\frac{n+2}2} - 2^{\frac{n+6}4} &\text{if $n\equiv2\Mod4$,}
\end{cases}$
\itemit{ii} $\fix_{\Gmn}(\si) = 
\begin{cases}
\ds \binom{\nt}{\mmt}-2\binom{\fnf}{\mmt} &\text{if $m\geq3$ is even}\\[5truemm]
\ds 2\binom{\nt-1}{\fmt}-2\binom{\fnf}{\fmt} &\text{if $m\geq3$ is odd.}
\end{cases}$
\eit
\end{lemma}

\pf
Write $n=2k$, and let $l=\lfloor\frac k2\rfloor$.  This time we may assume the fixed points of $\si$ are $1$ and $k+1$, in which case $\si(x)\equiv n+2-x\Mod n$ for all $x$ (cf.~Figure \ref{fig:circles4}).  Define a map $f:T_{k+1}\to T_n$ by $f(a_1,\ldots,a_{k+1})=(a_1,a_2,\ldots,a_k,a_{k+1},a_k,\ldots,a_2)$.

\pfitem{i}
As usual, $f$ is injective and has image $\Fix_{T_n}(\si)$, so that $\fix_{T_n}(\si)=2^{k+1}=2^{\frac{n+2}2}$.  Thus, by \eqref{eq:fix}, it remains to show that
\begin{equation}\label{eq:poly5}
\fix_{B_n}(\si) = \begin{cases}
\ds2^{\frac{n+8}4}-1 &\text{if $k=2l$}\\
\ds2^{\frac{n+6}4} &\text{if $k=2l+1$.}
\end{cases}
\end{equation}
To do so, consider some $n$-tuple $f(\ba)\in\Fix_{T_n}(\si)$, where $\ba\in T_{k+1}$.  Then $f(\ba)\in B_n$ if and only if at least one of the following holds (cf.~Figure \ref{fig:circles4}):
\[
\text{(a) $(a_1,\ldots,a_{k-l})=(0,\ldots,0)$, or \qquad\quad (b) $(a_{l+2},\ldots,a_{k+1})=(0,\ldots,0)$, or \qquad\quad (c) $\ba=(1,0,\ldots,0,1)$.}
\]
For $k=2l$ or $2l+1$, respectively, there are $2\cdot2^{l+1}-(2\text{ or }1)$ elements $\ba\in T_{k+1}$ satisfying (a) or (b).  Since $\ba=(1,0,\ldots,0,1)$ satisfies neither (a) nor (b), it follows that $\fix_{B_n}(\si) = 2^{l+1}-(2\text{ or }1)+1$, which reduces to \eqref{eq:poly5}.

\pfitem{ii} Write $h=\fmt$, and consider some $n$-tuple $f(\ba)\in\Fix_{\Tmn}(\si)$, where $\ba\in T_{k+1}$.  If $m=2h+1$ is odd, then exactly one of $a_1,a_{k+1}$ must be 1, and $h$ of the remaining entries $a_2,\ldots,a_k$ must be 1's.  If $m=2h$ is even, then either $a_1=a_{k+1}=1$ and $h-1$ of $a_2,\ldots,a_k$ are 1's, or else $a_1=a_{k+1}=0$ and $h$ of $a_2,\ldots,a_k$ are~1's.  Thus,
\[
\fix_{\Tmn}(\si) = \begin{cases}
\binom{k-1}{h-1}+\binom{k-1}h=\binom kh &\text{if $m=2h$ is even}\\
2\binom{k-1}h &\text{if $m=2h+1$ is odd.}
\end{cases}
\]
By \eqref{eq:fix}, and since $k=\frac n2$ and $h=\fmt$, it remains to show that $\fix_{\Bmn}(\si)=2\binom{\fnf}{\fmt}$.  To do so, consider an~$n$-tuple $f(\ba)\in\Fix_{\Bmn}(\si)$, where $\ba\in T_{k+1}$.  So $\ba$ satisfies either (a) or (b), as above, and these are mutually exclusive (since $m\geq3$); note that (c) cannot hold (also since $m\geq3$).  Note that $\ba$ satisfies (a) if and only if $a_{k+1}=0$ or $1$ (for $m=2h$ or $2h+1$, respectively), and $h$ of $a_{k-l+1},\ldots,a_k$ are 1's.  There are thus $\binom lh = \binom{\fnf}{\fmt}$ such $(k+1)$-tuples $\ba$ satisfying (a), and there are the same number satisfying (b).  
\epf

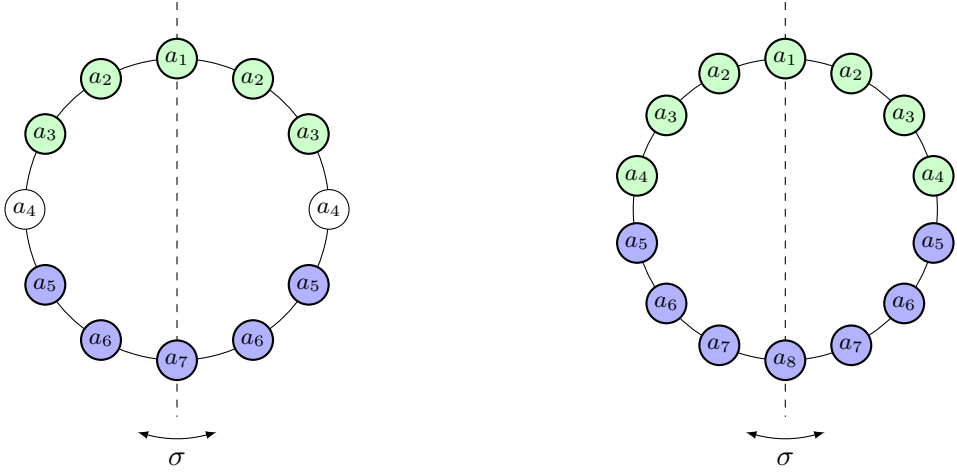
\begin{figure}[t]%
\begin{center}
\begin{tikzpicture}[scale=1]
\begin{scope}
\evencirclediagramnewnew{12}{11/3,12/2,1/1,2/2,3/3}{5/5,6/6,7/7,8/6,9/5}{4/4,10/4}
\end{scope}
\begin{scope}[shift={(8,0)}]
\evencirclediagramnewnew{14}{12/4,13/3,14/2,1/1,2/2,3/3,4/4}{5/5,6/6,7/7,8/8,9/7,10/6,11/5}{}
\end{scope}
\end{tikzpicture}
\end{center}
\vspace{-0.6cm}
\caption{The reflection $\si\in\D_n$ and the $n$-tuple $f(\ba)=(a_1,a_2,\ldots,a_k,a_{k+1},a_k,\ldots,a_2)\in\Fix_{T_n}(\si)$, from the proof of Lemma~\ref{lem:poly5}, in the cases $n=12$ (left) and $n=14$ (right).  If condition (a) or (b) from the proof is satisfied, then the green (light) or blue (dark) vertices, respectively, yield a bad block of 0's in $f(\ba)$.  If neither (a) nor~(b) holds, then~$f(\ba)$ has no bad blocks.}
\label{fig:circles4}
\end{figure}

\section{Proofs of the main results}\label{sect:proofs}

We are now in a position to complete the proofs of our main results, as stated in Section \ref{sect:intro}.  

\pf[\bf Proof of Theorem \ref{thm:pmn}]
As discussed in Section \ref{sect:dihedral}, the proof uses Burnside's Lemma.  Specifically, equation~\eqref{eq:BS} says that $p_{m,n} = \frac1{2n}\sum_{\si\in\D_n}\fix_{\Gmn}(\si)$.  From Lemmas \ref{lem:poly1} and \ref{lem:poly2}, and the fact that the cyclic group $\C_n$ has $\varphi(d)$ elements of order $d$ if $d\pre n$, the contribution to \eqref{eq:BS} made by the rotations from $\D_n$ is
\[
\frac1{2n}\sum_{\si\in\C_n}\fix_{\Gmn}(\si) = \frac1{2n} \left(
\binom nm - n\binom{\fnt}{m-1}
+ \sum_{1\not=d\pre n} \varphi(d)\binom{\frac nd}{\frac md}
\right)
= \sum_{d\pre n} \frac{\varphi(d)}{2n}\binom{\frac nd}{\frac md} - \frac12\binom{\fnt}{m-1}.
\]
Since $\binom{\frac nd}{\frac md}=0$ if $d\nmid m$, ``$\sum_{d\pre n}$'' may be replaced by ``$\sum_{d\pre\gcd(m,n)}$''.  It therefore remains to show that the contribution to \eqref{eq:BS} made by the reflections from $\D_n$ is
\[
\frac1{2n}\sum_{\si\in\D_n\sm\C_n}\fix_{\Gmn}(\si) = 
\frac12 \left(
\binom{\fmt+\lfloor\frac{n-m}2\rfloor}{\fmt}
- \binom{\fnf}{\fmt}
- \binom{\ntf}{\mmt}
\right).
\]
If $n$ is odd, then $\D_n$ has $n$ reflections; if $n$ is even, then $\D_n$ has $\nt$ reflections with no fixed points, and $\nt$ with two fixed points.   Combining Lemmas \ref{lem:poly3}, \ref{lem:poly4} and \ref{lem:poly5}, and keeping in mind that $\binom{\ntf}{\mmt}=0$ if $m$ is odd, we calculate case-by-case that
\[
\frac1{2n}\sum_{\si\in\D_n\sm\C_n}\fix_{\Gmn}(\si) = 
\begin{cases}
\ds \frac12 \left(\binom{\fnt}{\fmt} - \binom{\fnf}{\fmt} - \binom{\ntf}{\frac m2}\right) &\text{if $n$ is odd}\\[6truemm]
\ds \frac12 \left(\binom{\nt}{\mmt} - \binom{\fnf}{\mmt} - \binom{\ntf}{\mmt}\right) &\text{if $n$ and $m$ are both even}\\[6truemm]
\ds \frac12 \left(\binom{\nt-1}{\fmt} - \binom{\fnf}{\fmt} - \binom{\ntf}{\mmt}\right) &\text{if $n$ is even $m$ is odd.}
\end{cases}
\]
It finally remains to observe that $\fmt+\lfloor\frac{n-m}2\rfloor$ is equal to $\frac n2-1$ if $n$ is even and $m$ is odd, or to $\fnt$ otherwise; indeed, this may be checked on a case-by-case basis.
\epf

\pf[\bf Proof of Theorem \ref{thm:poly}]
With some careful algebra, \eqref{eq:BS} and Lemmas \ref{lem:poly1}--\ref{lem:poly5} yield
\begin{equation}\label{eq:poly_complicated}
p_n = \sum_{d\pre n}\frac{\varphi(d)}{2n}\big(2^{\frac nd}-1\big) - 2^{\lfloor\frac n2\rfloor-1} +
\begin{cases}
3\cdot2^{\frac{n-4}2}-3\cdot2^{\frac{n-4}4}+\frac12 &\text{if $n\equiv0\Mod4$}\\[3truemm]
2^{\frac{n-1}2}-3\cdot2^{\frac{n-5}4}+\frac12 &\text{if $n\equiv1\Mod4$}\\[3truemm]
3\cdot2^{\frac{n-4}2}-2^{\frac{n+2}4}+\frac12 &\text{if $n\equiv2\Mod4$}\\[3truemm]
2^{\frac{n-1}2}-2^{\frac{n+1}4}+\frac12 &\text{if $n\equiv3\Mod4$.}
\end{cases}
\end{equation}
(Note that the ``$\nt$'' from the even case of Lemma \ref{lem:poly1}(i) cancels out with the ``$\nt$'' from the $d=2$ case of Lemma \ref{lem:poly2}(i).)
Using the identity $\sum_{d\pre n}\varphi(d)=n$, we have
\[
\sum_{d\pre n}\frac{\varphi(d)}{2n}\big(2^{\frac nd}-1\big) = \sum_{d\pre n}\frac{\varphi(d)\cdot2^{\frac nd}}{2n} - \frac1{2n}\sum_{d\pre n}\varphi(d) = \sum_{d\pre n}\frac{\varphi(d)\cdot2^{\frac nd-1}}{n} - \frac12.
\]
The ``$\frac12$'' here cancels those in \eqref{eq:poly_complicated}.  If $n$ is even, then $- 2^{\lfloor\frac n2\rfloor-1} + 3\cdot2^{\frac{n-4}2} = - 2\cdot2^{\frac{n-4}2} + 3\cdot2^{\frac{n-4}2} = 2^{\frac{n-4}2} = 2^{\lfloor\frac{n-3}2\rfloor}$.  A similar calculation shows that if $n$ is odd, then $- 2^{\lfloor\frac n2\rfloor-1} + 2^{\frac{n-1}2} = 2^{\frac{n-3}2} = 2^{\lfloor\frac{n-3}2\rfloor}$.
\epf

\pf[\bf Proof of Theorem \ref{thm:asymptotic}]
The dominant term in Theorem \ref{thm:poly} is the $d=1$ term of the sum: i.e.,~$\frac{2^{n-1}}n$.  All other terms (of which there are at most $n$) are at most a constant multiple of~$2^{\frac n2}$.  
\epf

The formula $p_n\sim\frac{2^{n-1}}n=\frac{2^n}{2n}$ given in Theorem \ref{thm:asymptotic} can be interpreted as saying that: (i) practically all of the $2^n$ subsets of $\{1,\ldots,n\}$ correspond to polygons (cf.~Lemma \ref{lem:poly1}(i)); and (ii) practically all such polygons are completely asymmetric, so that each equivalence class of polygons is counted approximately~$2n$ times.

\pf[\bf Proof of Theorem \ref{thm:asymptotic2}]
Note that for integers $0\leq k\leq x$, the binomial coefficient $\binom xk$ is a polynomial of degree $k$ in $x$; thus, it is asymptotic to its leading term: $\binom xk\sim \frac{x^k}{k!}$.  Examining the expression for $p_{m,n}$ in Theorem \ref{thm:pmn}, there are two dominant terms: the $d=1$ term of the sum, and the second bracketed binomial coefficient, both of which are polynomials in $n$ of degree $m-1$.  Thus, for fixed $m$, and as $n\to\infty$,
\[
p_{m,n} 
\sim \frac1{2n}\binom nm - \frac12\binom{\fnt}{m-1} 
\sim \frac1{2n}\cdot\frac{n^m}{m!} - \frac12\cdot\frac{\frac{n^{m-1}}{2^{m-1}}}{(m-1)!}
=\frac{n^{m-1}}{2m!} - \frac{n^{m-1}}{2^m(m-1)!}
= \frac{2^{m-1} - m}{2^mm!}  n^{m-1}.
\qedhere
\]
\epf

\section{Triangles and quadrilaterals}\label{sect:triquad}

Specialising Theorem \ref{thm:pmn} to the cases $m=3$ and $m=4$ allows us to recover Honsberger's Theorem on triangles (stated in Theorem \ref{thm:tn} above), and also the new result (Theorem \ref{thm:qn}) on quadrilaterals.

\pf[\bf Proof of Theorem \ref{thm:tn}]
By Theorem \ref{thm:pmn}, the number of inequivalent (i.e., incongruent) integer triangles of perimeter $n$ is 
\begin{align*}
p_{3,n} &= \frac{1}{2n}\binom n3 + \frac{1}{n}\binom{\frac n3}1 +\frac12\left(
\binom{1+\lfloor\frac{n-3}2\rfloor}1 - \binom{\fnt}2 - \binom{\fnf}1-0
\right) \\[3truemm]
&= \frac{n^2-3n+2}{12} + \frac{(0\text{ or }1)}3 +
\frac{1+\lfloor\frac{n-3}2\rfloor}2 - \frac{\fnt(\fnt-1)}4 - \frac{\fnf}2
\\[3truemm]
&= \begin{cases}
\ds \frac{n^2}{48} + \frac{-4+(0\text{ or }4)+(0\text{ or }3)}{12} &\text{if $n$ is even}\\[5truemm]
\ds \frac{(n+3)^2}{48} + \frac{-4+(0\text{ or }4)+(0\text{ or }3)}{12} &\text{if $n$ is odd.}
\end{cases}
\end{align*}
Since $p_{3,n}$ is an integer, and since $-\frac12<-\frac{4}{12}\leq\frac{-4+(0\text{ or }4)+(0\text{ or }3)}{12}\leq\frac{3}{12}<\frac12$, it follows that $p_{3,n}$ is the nearest integer to $\frac{n^2}{48}$ or $\frac{(n+3)^2}{48}$, as appropriate.
\epf

\pf[{\bf Proof of Theorem \ref{thm:qn}}]
As in the above proof of Theorem \ref{thm:tn}, we use Theorem \ref{thm:pmn} to show, case by case, that 
\begin{equation}\label{eq:quad}
p_{4,n} = \begin{cases}
\frac{n^3-3n^2+20n}{96} &\text{if $n\equiv0\Mod4$}\\
\frac{n^3-7n+6}{96} &\text{if $n\equiv1\Mod4$}\\
\frac{n^3-3n^2+20n-36}{96} &\text{if $n\equiv2\Mod4$}\\
\frac{n^3-7n-6}{96} &\text{if $n\equiv3\Mod4$.}
\end{cases}
\end{equation}
For example, if $n\equiv2\Mod4$, then
{\small
\begin{align*}
p_{4,n} &= \frac1{2n}\binom n4 + \frac1{2n}\binom{\frac n2}2 + \frac12 \left(\binom{\frac n2}2 - \binom{\frac n2}3 - \binom{\frac{n-2}4}2 - \binom{\frac{n+2}4}2\right) \\[5truemm]
&= \frac1{2n}\cdot\frac{n(n-1)(n-2)(n-3)}{24} 
+ \frac1{2n}\cdot\frac{\frac n2\cdot\frac{n-2}2}2 + \frac12\left(
\frac{\frac n2\cdot\frac{n-2}2}2 
- \frac{\frac n2\cdot\frac{n-2}2\cdot\frac{n-4}2}6 
- \frac{\frac {n-2}4\cdot\frac{n-6}4}2
- \frac{\frac {n+2}4\cdot\frac{n-2}4}2 \right),
\end{align*}
}
\!\!which reduces to $\frac{n^3-3n^2+20n-36}{96}$.  Equation \eqref{eq:quad} clearly completes the proof in the $n\equiv0\Mod4$ case.  In the $n\equiv1$ case, since $\frac{n^3-7n+6}{96}$ is an integer, and since $\frac6{96}<\frac12$, certainly $\frac{n^3-7n+6}{96}$ is the nearest integer to~$\frac{n^3-7n}{96}$.  The other cases are analogous.
\epf

\section{Calculated values and concluding remarks}\label{sect:values}

Tables \ref{tab:values} and \ref{tab:sums} below give calculated values of $p_{m,n}$ and~$p_n$, respectively, for $3\leq n\leq20$.  Because of the powers of~2 in Theorem \ref{thm:poly}, the $p_n$ sequence is very easy to compute; for example, the millionth term can be calculated in well under 20 seconds on a standard laptop.
The first four rows of Table \ref{tab:values} are Sequences A005044, A057886, A124285 and A124286, respectively, on the Online Encyclopedia of Integer Sequences \cite{OEIS}, while the whole of Table~\ref{tab:values} is Sequence A124287.  At the time of writing, the sequence~$p_n$ (Table \ref{tab:sums}) did not appear on \cite{OEIS}; it is now Sequence A293818.  Out of all these sequences, as far as we are aware, a formula had only previously been proven for Sequence A005044 (Honsberger's Theorem~\ref{thm:tn} concerning triangles, quoted in Section~\ref{sect:intro} above).

If we were interested in polygon equivalence under cyclic re-ordering of edges but not reversals, then the inequivalent integer $m$-gons (respectively, polygons) of perimeter $n$ are in one-one correspondence with the orbits of $\Gmn$ (respectively, $G_n$) under the action of the cyclic group $\C_n$ given by \eqref{eq:action}.  If we write $p'_{m,n}$ (respectively, $p'_n$) for the number of such inequivalent $m$-gons (respectively, polygons), then
\begin{align*}
p_{m,n}' = |\Gmn/\C_n| &= \frac1n \sum_{\si\in\C_n}\fix_{\Gmn}(\si) = \sum_{d\pre \gcd(m,n)}\frac{\varphi(d)}{n}\binom{\frac nd}{\frac md} -\binom{\fnt}{m-1} , \\[5truemm]
\text{and}\qquad p_n' = |G_n/\C_n| &= \frac1n \sum_{\si\in\C_n}\fix_{G_n}(\si) = \sum_{d\pre n}\frac{\varphi(d)}{n}\big(2^{\frac nd}-1\big) - 2^{\lfloor\frac n2\rfloor} = \sum_{d\pre n}\frac{\varphi(d)\cdot2^{\frac nd}}{n}-1 - 2^{\lfloor\frac n2\rfloor}.
\end{align*}
The sequence $p_{3,n}'$ (triangles up to rotation) is Sequence A008742 on \cite{OEIS}; curiously, however, A008742 is listed as the Molien series of a certain three dimensional point group.  More computed values may be found in Sequences A293819--A293823 on \cite{OEIS}.
(Note that Sequence~A124278 of \cite{OEIS} counts $|\Pmn/\S_m|$: i.e., polygons up to arbitrary reorderings of the sides, under which $(1,1,2,2)$ and $(1,2,1,2)$ are considered equivalent, for example; see also \cite{APR2001}.)

Finally, we observe that Theorem \ref{thm:poly} has an interesting number-theoretic consequence: namely, that $\sum_{d\pre n}\frac{\varphi(d)\cdot2^{\frac nd-1}}n$ is an integer for $n\geq3$; cf.~Sequences A000031, A053634 and A053635 on \cite{OEIS}.

\begin{table}[h]
\begin{center}
\begin{tabular}{rrrrrrrrrrrrrrrrrrr}
\hline
$m\sm n$ & 3 & 4 & 5 & 6 & 7 & 8 & 9 & 10 & 11 & 12 & 13 & 14 & 15 & 16 & 17 & 18 & 19 & 20 \\
\hline
3$\phantom{{}\sm n}$ & 1 & 0 & 1 & 1 & 2 & 1 & 3 & 2 & 4 & 3 & 5 & 4 & 7 & 5 & 8 & 7 & 10 & 8 \\
4$\phantom{{}\sm n}$ & & 1 & 1 & 2 & 3 & 5 & 7 & 9 & 13 & 16 & 22 & 25 & 34 & 38 & 50 & 54 & 70 & 75 \\
5$\phantom{{}\sm n}$ & & & 1 & 1 & 3 & 4 & 9 & 13 & 23 & 29 & 48 & 60 & 92 & 109 & 158 & 186 & 258 & 296 \\
6$\phantom{{}\sm n}$ & & & & 1 & 1 & 4 & 7 & 15 & 25 & 46 & 72 & 113 & 172 & 248 & 360 & 491 & 686 & 896 \\
7$\phantom{{}\sm n}$ & & & & & 1 & 1 & 4 & 8 & 20 & 37 & 75 & 129 & 228 & 359 & 584 & 868 & 1324 & 1870 \\
8$\phantom{{}\sm n}$ & & & & & & 1 & 1 & 5 & 10 & 29 & 57 & 125 & 231 & 435 & 745 & 1261 & 2031 & 3195 \\
9$\phantom{{}\sm n}$ & & & & & & & 1 & 1 & 5 & 12 & 35 & 79 & 185 & 374 & 749 & 1382 & 2489 & 4237 \\
10$\phantom{{}\sm n}$ & & & & & & & & 1 & 1 & 6 & 14 & 47 & 111 & 280 & 600 & 1281 & 2493 & 4746 \\
11$\phantom{{}\sm n}$ & & & & & & & & & 1 & 1 & 6 & 16 & 56 & 147 & 392 & 912 & 2052 & 4261 \\
12$\phantom{{}\sm n}$ & & & & & & & & & & 1 & 1 & 7 & 19 & 72 & 196 & 561 & 1368 & 3260 \\
13$\phantom{{}\sm n}$ & & & & & & & & & & & 1 & 1 & 7 & 21 & 84 & 252 & 756 & 1980 \\
14$\phantom{{}\sm n}$ & & & & & & & & & & & & 1 & 1 & 8 & 24 & 104 & 324 & 1032 \\
15$\phantom{{}\sm n}$ & & & & & & & & & & & & & 1 & 1 & 8 & 27 & 120 & 406 \\
16$\phantom{{}\sm n}$ & & & & & & & & & & & & & & 1 & 1 & 9 & 30 & 145 \\
17$\phantom{{}\sm n}$ & & & & & & & & & & & & & & & 1 & 1 & 9 & 33 \\
18$\phantom{{}\sm n}$ & & & & & & & & & & & & & & & & 1 & 1 & 10 \\
19$\phantom{{}\sm n}$ & & & & & & & & & & & & & & & & & 1 & 1 \\
20$\phantom{{}\sm n}$ & & & & & & & & & & & & & & & & & & 1 \\
\hline
\end{tabular}
\end{center}
\vspace{-5truemm}
\caption{The number $p_{m,n}$ of inequivalent integer $m$-gons with perimeter $n$, for $3\leq m\leq n\leq20$.}
\label{tab:values}
\end{table}

\begin{table}[h]
\begin{center}
\begin{tabular}{rrrrrrrrrrrrrrrrrrr}
\hline
$n$ & 3 & 4 & 5 & 6 & 7 & 8 & 9 & 10 & 11 & 12 & 13 & 14 & 15 & 16 & 17 & 18 & 19 & 20 \\
\hline
& 1 & 1 & 3 & 5 & 10 & 16 & 32 & 54 & 102 & 180 & 336 & 607 & 1144 & 2098 & 3960 & 7397 & 14022 & 26452 \\
\hline
\end{tabular}
\end{center}
\vspace{-5truemm}
\caption{The number $p_n$ of inequivalent integer polygons with perimeter $n$, for $3\leq n\leq20$.  These are column sums of Table \ref{tab:values}.}
\label{tab:sums}
\end{table}

\footnotesize
\def\bibspacing{-1.1pt}
\bibliography{biblio}
\bibliographystyle{plain}

~

\noindent JE: Centre for Research in Mathematics, Western Sydney University, Sydney, Australia; {\tt J.East\,@\,WesternSydney.edu.au}
\\[2mm]
\noindent RN: {\tt rniles\,@\,yahoo.com}

\end{document}